\documentclass[12pt,oneside]{article}
\usepackage{amsmath,amssymb,amsfonts,amsthm}
\usepackage{color}
\usepackage[all]{xy}
 \usepackage{color}
\usepackage{makeidx}
\usepackage{multirow}
\usepackage{rotating}
\usepackage{pdfpages}
\usepackage{pdflscape}
\usepackage{tabularray}
\usepackage{multicol}
\usepackage{amscd}

\allowdisplaybreaks

\newcommand{\set}[1]{\left\{#1\right\}}

\def\paa{\dot{\partial}}

\numberwithin{equation}{section} 
\numberwithin{figure}{section} 

\theoremstyle{plain}
\newtheorem*{theorem*}{Theorem}
\newtheorem{theorem}{Theorem}[section]
\newtheorem{lemma}[theorem]{Lemma}
\newtheorem{proposition}[theorem]{Proposition}
\newtheorem{corollary}[theorem]{Corollary}
\newtheorem{remark}[theorem]{Remark}

\newtheorem{definition}[theorem]{Definition}

\newtheorem{example}{Example}
\newtheorem*{acknowledgement*}{Acknowledgement}

\numberwithin{equation}{section}

\newcommand\overcirc[1]{\raisebox{10pt}{\tiny{$\circ$}}{\kern-7.5pt}\mbox{$#1$}}
\newcommand\undersym[2]{\raisebox{-6pt}{$#2$}{\kern-5pt}\mbox{$#1$}}
\newcommand\overdiamond[1]{\raisebox{10pt}{\small$\star$}{\kern-7.5pt}\mbox{$#1$}}
\newcommand\overast[1]{\raisebox{10pt}{\small$\ast$}{\kern-7.5pt}\mbox{$#1$}}
\newcommand\overlind[1]{\raisebox{10pt}{\small$\overline{{\hspace{2pt}}\star}$}{\kern-7.5pt}\mbox{$#1$}}
\newcommand\overlinc[1]{\raisebox{10pt}{\small$\overline{{\hspace{2pt}}\circ}$}{\kern-7.5pt}\mbox{$#1$}}
\newcommand\overlina[1]{\raisebox{10pt}{\small$\overline{{\hspace{1pt}}\ast}$}{\kern-7.5pt}\mbox{$#1$}}

\newcommand\undersymm[2]{\raisebox{-7pt}{\tiny$#2$}{\kern-15pt}\mbox{$#1$}}

\begin{document}

\title{\textbf{On the covariant coefficients of geodesic  sprays on  Finsler manifolds} }
\author{\bf{S. G. Elgendi$^{\,1}$, A. Soleiman$^{\,2}$ and  Nabil L. Youssef$^{\,3}$}}
\date{}
\maketitle                     
\vspace{-1.16cm}

\begin{center}
{ $^{1}$ Department of Mathematics, Faculty of Science, Islamic University of Madinah,  Madinah, Kingdom of Saudia Arabia
\vspace{0.2cm}\\
$^{2}$  Department of  Mathematics, Faculty of Science, Al Jouf University, Skaka, Kingdom of Saudia Arabia
\vspace{0.2cm}\\
$^{3}$ Department of Mathematics, Faculty of Science, Cairo
University, Giza,  {Egypt}}
\end{center}

\begin{center}
E-mails: selgendi@iu.edu.sa, salahelgendi@yahoo.com\\
asoliman@ju.edu.sa, amrsoleiman@yahoo.com\\
nlyoussef@sci.cu.edu.eg, nlyoussef2003@yahoo.fr  
\end{center} 

\maketitle                     
\vspace{-1.0cm} 

\vspace{0.6cm}
\begin{abstract} 
For a Finsler metric $F$, we introduce the notion  of $F$-covariant coefficients $H_i$ of the geodesic spray of $F$ (Def. 3.1). We study some geometric consequences concerning the objects $H_i$. If the $F$-covariant coefficients  $H_i$ are written in the form $H_i={\dot{\partial}}_iH$, for some smooth function $H$ on ${\mathcal T\hspace{-1pt}M}$, positively 3-homogeneous in y, then $H$ is called spray scalar or simply $S$-scalar. We prove that if the $S$-scalar exists, then it is of the form $H=\frac{1}{12}\,y^i\partial_iF^2$ and this expression is unique up to a function of position only. We prove also that on a Finsler maifold $(M,F)$, the $S$-scalar $H$ exists if and only if $(M,F)$ is dually flat. 
Generally, the $n^3$ functions  $H^h_{ij}$ resulting from the $F$-covariant coefficients do not  form a linear connection. We find out that in the case of projectively flat metrics, the $n^3$ functions $H^h_{ij}$ are coefficients of a linear connection. 
We introduce two new special Finsler spaces, namely, the $H$-Berwald and the $H$-Landsberg spaces and show that every $H$-Berwald metric is $H$-Landsbergian but the converse is not necessarily true. 
Also, we study the $F$-covariant coefficients $H_i$ of projectivly  flat and dually flat spherically symmetric Finsler metrics and provide a solution of the \lq\lq$H$-unicorn" Landsberg problem.
Finally, we give some examples of $H$-Berwald and $H$-Landsberg metrics and an example of $H$-Landsberg metric which is not $H$-Berwaldian.    
\end{abstract}
 



\section{Introduction}

~\par

W. Ambrose et al. introduced the concept of sprays in 1960 \cite{sprays}. A system of second-order ordinary differential equations (SODEs) with positively 2-homogeneous coefficient functions can be interpreted as a second-order vector field, known as a spray. Specifically, for an $n$-dimensional differentiable manifold $M$, a spray on $M $ is a vector field on the tangent bundle $T M$. The geodesic spray of a Riemannian or Finslerian metric $F$ corresponds to the geodesic equation of that metric. For more details on the geometry of sprays, refer to \cite{r21, Szilasi-book}.

\medskip

In this paper, we start with a Finsler manifold $(M,F)$ and its geodesic spray $S$.  Then, we introduce the notion of $F$-covariant coefficients $H_i$ of $S$ as follows
$$
 H_{i}= \frac{1}{4}(y^r\partial_r\dot{\partial}_iF^2 - \partial_iF^2).
$$
 We study some geometric consequences related to  the objects $H_i$.  Moreover, if $H_i$ are obtained from a scalar function $H$ on $TM$ by differentiating    with respect to the fiber variables $y^i$, then we call $H$ an $S$-scalar. For a  Finsler manifold  $(M,F)$    with   geodesic spray   $S$, we prove that  the $S$-scalar $H$, if exists,  is given by   
$$
 H(x, y)=\frac{1}{12} y^k \partial_k F^2.
$$
Moreover, the existence of the $S$-scalar $H$ is characterized by the  dual flatness of $(M,F)$.

\medskip

By taking derivatives of $H_i$ with respect to fiber variables,  we get $H_{ij}$, $H_{ijk}$ and $H_{ijkh}$, then we have the $n^3$ functions $H^i_{jk}:=H_{rjk} \ g^{ri}$. We find that   the $n^3$ functions $H^i_{jk}$ are not the coefficients of a linear connection, but   in the case of projectively flat metrics  they are. Moreover, we introduce  two new special Finsler spaces: 
 $H$-Berwald space  (the $H$-Berwald tensor $H_{ijkh}$ vanishes) and
  $H$-Landsberg  space (the $H$-Landsberg tensor $\mathcal{L}_{jkh}:=y^iH_{ijkh}$ vanishes). We  show that every $H$-Berwald metric is $H$-Landsbergian but the converse is not necessarily true.
 
  We investigate the $F$-covariant coefficients $H_i$ of a spherically symmetric Finsler metric and then focus our attention  on the class of projectively flat and dually flat spherically symmetric Finsler metrics characterized in \cite{Najafi}. This class represents  a solution of the \lq\lq $H$-unicorn" Landsberg problem.  We provide some examples of  $H$-Berwald and  $H$-Landsberg metrics and, furthermore, we give an example of $H$-Landsberg metric which is not $H$-Berwaldian.    
 Namely, we consider the Randers metric $F$ on $\mathbb{R}^n$ given by
 $$F(x,y)=c\sqrt{y_1^2+y_2^2+\cdots + y_n^2}+(a_1+x_1)y^1+(a_2+x_2)y^2+\cdots+(a_n+x_n)y^n,$$
where  $c\in\mathbb{R}$ and $a$ is a fixed vector in $\mathbb{R}^n$.
The metric  $F$  is $H$-Landsbergian metric but not $H$-Berwaldian.    

It should be noted that, whenever suitable, we use the  Finsler package \cite{NF_Package} for the calculations of the examples throughout.

   \section{Preliminaries }

Let $M$ be an $n$-dimensional smooth manifold with tangent bundle $\left(T M, \pi_M, M\right)$ and subbundle of nonzero tangent vectors $(\mathcal{T} M, \pi, M)$. Consider local coordinates $\left(x^i\right)$ on $M$ and corresponding coordinates $\left(x^i, y^i\right)$ on $T M$. The natural almost-tangent structure $J$ on $T M$ is a vector 1 -form locally given by $J=\frac{\partial}{\partial y^i} \otimes d x^i$. The canonical or Liouville vector field $\mathcal{C}$ on $T M$ is the vertical vector field defined by $\mathcal{C}=y^i \frac{\partial}{\partial y^i}$.

A vector field $S\in \mathfrak{X}(\mathcal{T} M)$ is    a spray if $JS = \mathcal{C}$ and $[\mathcal{C}, S] = S$.  A spray is  written, locally, in the form
\begin{equation}
  \label{eq:spray}
  S = y^i \frac{\partial}{\partial x^i} - 2G^i\frac{\partial}{\partial y^i},
\end{equation}
where the functions $G^i=G^i(x, y)$ are   the spray coefficients, which are positively 2-homogeneous in $y$.

A nonlinear connection is    an $n$-dimensional distribution $H_z(\mathcal{T} M)$ that complements the vertical distribution $V_z(\mathcal{T} M)$ for all $z \in \mathcal{T} M$. Therefore, we have: 
$$T_z(\mathcal{T}M) = H_z(\mathcal{T}M) \oplus V_z(\mathcal{T}M).$$

Every spray $S$ defines a canonical nonlinear connection through its associated horizontal and vertical projectors
\begin{equation}
  \label{projectors}
    h=\frac{1}{2} (id_{TM} + [J,S]), \,\,\,\,            v=\frac{1}{2}(id_{TM} - [J,S])
\end{equation}

 Additionally, the canonical nonlinear connection attached to a spray can be characterized in terms of an almost product structure $\Gamma=$ $[J, S]=h-v$.  Locally, the two projectors $h$ and $v$ are expressed as
$$h=\frac{\delta}{\delta x^i}\otimes dx^i, \quad\quad v=\frac{\partial}{\partial y^i}\otimes \delta y^i.$$
Moreover, 
$$\frac{\delta}{\delta x^i}:=\frac{\partial}{\partial x^i}-N^j_i(x,y)\frac{\partial}{\partial y^j},\quad \delta y^i:=dy^i+N^j_i(x,y)dx^i, \quad N^j_i(x,y):=\frac{\partial G^j}{\partial y^i}.$$
 The Berwald connection is a fundamental object in Finsler geometry, and it plays a similar role to the Levi-Civita connection in Riemannian geometry. The coefficients $G^h_{ij}$ of  the  Berwald connection are determined  by the formula
$$ G^h_{ij}=\frac{\partial G^h_j}{\partial y^i}.$$
For simplicity,  we use the notations
$$\partial_j:=\frac{\partial}{\partial x^j}, \quad \dot{\partial}_j:=\frac{\partial}{\partial y^j}.$$

 A Finsler structure on $M$ is defined as follows.

\begin{definition}\label{fin.struc.} A Finsler structure, or a Finsler metric, on a smooth manifold $M$ is a function 
$$F:TM\rightarrow \mathbb{R}$$
satisfying the conditions:
\begin{description}
   \item[(a)] $F\geqslant 0$ and $F(x,y)=0$ if and only if $y=0$.

    \item[(b)] $F$ is $C^\infty$ on the slit tangent
    bundle  $\mathcal{T}M=TM\backslash\{0\}$.

    \item[(c)] $F(x,y)$ is positively homogeneous of degree one in $y$: $F(x,\lambda y) = \lambda F(x,y)$
     for all $y \in TM$ and $\lambda > 0$.

    \item[(d)] The Hessian matrix ${g_{ij}(x,y):=\frac{1}{2}\dot{\partial}_i\dot{\partial}_j F^2}$
is non-degenerate at each point of $\mathcal{T}M$.
\end{description}
The pair $(M,F)$ is called  a Finsler manifold, and the tensor  $g=g_{ij}(x,y)\,dx^i\otimes dx^j$ is called the
Finsler metric  tensor  of   $(M,F)$.
\end{definition}

As the $2$-form $dd_JE$; $E:=\frac{1}{2}F^2$, is
non-degenerate,   Grifone \cite{r21}  proved the following fundamental result.

\begin{proposition}\cite{r21}
  The  Euler-Lagrange equation
\begin{equation*}
  \label{eq:EL}
   i_Sdd_JE+dE =0
\end{equation*}
This determines a unique spray $S$ on $M$, known as the geodesic spray of the Finsler structure $F$.  
\end{proposition}

\begin{definition}
 A spray $S$ on a manifold $M$ is   Finsler metrizable if there exists a Finsler metric $F$ such that the geodesics of $F$ coincide with the integral curves of $S$.  In other words, the spray $S$ can be realized as the geodesic spray of a Finsler metric.
\end{definition}

For a Finsler manifold   $(M,F)$, the coefficients $G^i$ of the geodesic spray  of $F$ are given  by
 \begin{equation}
 \label{Eq:G^i}
 G^i= \frac{1}{4}g^{ih}( y^r\partial_r\dot{\partial}_hF^2 - \partial_hF^2).
 \end{equation}

\section{Covariant coefficients of a spray}

 For a Finsler manifold $(M,F)$ with   geodesic spray $S$, we focus our attention  to    the geometry of the objects   obtained from the the   coefficients $G^i$ of $S$ by lowering the index $i$. For this purpose, we set the following definition.

\begin{definition}
Let $(M,F)$ be a Finsler manifold with   geodesic spray $S$. The $F$-covariant coefficients $H_i$ of $S$ are defined by
\begin{equation}
\label{Eq:H_i}
 H_{i}:= g_{ir} G^r=\frac{1}{4}(y^r\partial_r\dot{\partial}_iF^2 - \partial_iF^2).
 \end{equation}
\end{definition}

Starting from the coefficients $G^i$ \eqref{Eq:G^i} of a metrizable spray $S$, we have a lot of associated geometric objects, for example, non-linear connection, linear connections and so on. Now, the question is that what kind of objects that can be obtained or constructed from the $F$-covariant coefficients $H_i$ of $S$.

\begin{lemma}
Under a change of coordinates $(x^i)\to (\widetilde{x}^i)$, the coefficients $H_i$ have  the transformation law
 \begin{equation}
\widetilde{H}_i=\frac{\partial {x}^h}{\partial \widetilde{x}^i}H_h-\frac{1}{2}\frac{\partial^2 \widetilde{x}^r}{\partial {x}^h\partial {x}^j} {y}^h {y}^j \widetilde{g}_{ir}
\end{equation}
\end{lemma}

It should be noted that  we denote $g_{ir}G^r$ by $H_i$ and not $G_i$ because when the spray $S$ is metrizable by essentially different Finsler metrics, then the corresponding   covariant coefficients may be different. To clarify this point, consider the following example.

\begin{example}
Let $M=\mathbb{B}^2\subset\mathbb{R}^2$, and  $\displaystyle{z_i=\frac{(1+\langle a,x\rangle)y_i-\langle a,y \rangle x_i}{\langle a,y \rangle}}$, where  $y\in T_x\mathbb{B}^2=\mathbb{R}^2$, $a=(a_1,a_2)\in \mathbb{R}^2$ is a constant vector with $|a|<1$.

Consider the Finsler metrics  $$F=\frac{\langle a,y\rangle \sqrt{z_1^2+z_2^2}}{(1+\langle a,x\rangle)^2}, \quad \overline{ F}= \frac{\langle a,y\rangle \sqrt{1+z_1^2+z_2^2}}{(1+\langle a,x\rangle)^2}, $$
 The  spray coefficients are given by
$$
G^i=-\frac{\langle a,y \rangle}{1+\langle a,x \rangle}y^i, \quad \overline{G}^i=G^i=-\frac{\langle a,y \rangle}{1+\langle a,x \rangle}y^i.
$$
For simplicity, let $a_1=0$ and $a_2=a$. The $F$-covariant coefficients $H_i$ are given by
$$H_1=\frac{4a y_1(a^2 x_1x_2y_2-a^2x_2^2y_1+ax_2y_2-y_1)}{(ax_1+1)^5}, \quad H_2=-\frac{4a y_1(a x_1y_2-ax_2y_1+y_2)}{(ax_1+1)^4},$$
 whereas the  $\overline{ F}$-covariant coefficients $\overline{ H}_i$ are given by
$$\overline{ H}_1=\frac{4a y_1(a^2 x_1x_2y_2-a^2x_2^2y_1-a^2y_1+ax_2y_2-y_1)}{(ax_1+1)^5}, \quad \overline{ H}_2=-\frac{4a y_1(a x_1y_2-ax_2y_1+y_2)}{(ax_1+1)^4}.$$
\end{example}

\begin{definition}
 Let $(M,F)$ be a   Finsler manifold with   geodesic spray $S$. If the   $F$-covariant coefficients $H_i$ are written in the form
 \begin{equation}
 \label{Eq:paa_iH}
 H_i=\dot{\partial}_i H, 
 \end{equation}
 for some $C^\infty$ scalar function  $H$ on $\mathcal{T}M$, positively 3-homogeneous in y, then $H$  is called  spray scalar or simply  $S$-scalar. 
\end{definition}

\begin{definition}\cite{Shen-dually}
A Finsler metric $F=F(x, y)$ on a smooth manifold $M$ is called  locally dually flat if at each point, there is a local coordinate system $\left(x^i\right)$ in which $F$ satisfies 
\begin{equation}\label{eq:Dually_flat_1}
y^k\dot{\partial}_i \partial_k F^2-2\partial_i F^2=0.
\end{equation}
Such a local system is called an adapted local system.\\
The Finsler metric  $F$ is globally dually flat, or simply dually flat, if it is locally dually flat and, in any adapted local system, satisfies the condition
\begin{equation}\label{eq:Dually_flat_2}
G^i=-\frac{1}{2} g^{i j} \dot{\partial}_j K,
\end{equation}
where $K=K(x, y)$ is a smooth function on $\mathcal{T} M$, positively $3$-homogeneous  in $y$ .
\end{definition}

\begin{proposition}\label{Prop:S-scalar} If the S-scalar $H$ exists, then it is of the form 
\begin{equation}\label{Eq:spray_scalar}
H=\frac{1}{12}\, y^i\partial_iF^2.
\end{equation}
The expression of $H$ in \eqref{Eq:spray_scalar} is unique up to a function of position only
\end{proposition}
\begin{proof}
Assume that the $S$-scalar $H$ exists, then $H_i=\dot{\partial}_iH$, by \eqref{Eq:paa_iH}. Contracting both sides of Equation \eqref{Eq:H_i} by $y^i$, we get
$$3H=y^i\dot{\partial}_iH=y^iH_i=\frac{1}{4} (y^r\partial_r(y^i\dot{\partial}_iF^2)-y^i\partial_iF^2)=\frac{1}{4}(2y^r\partial_rF^2-y^i\partial_iF^2)=\frac{1}{4}y^i\partial_iF^2.$$
Hence, $H=\frac{1}{12\,}y^i\partial_iF^2$.\\
For the uniqueness of the $S$-scalar $H$, let $H_1$ be another $S$-scalar. Then, $H_i=\dot{\partial}_iH_1$. Consequently, $\dot{\partial}_i(H_1-H)=\dot{\partial}_iH_1-\dot{\partial}_iH=0$. It follows that $H_1-H$ is a function of position only, say $f(x)$. Therefore, $H_1=H+f(x)=\frac{1}{12}\,y^i\partial_iF^2+f(x)$. 
\end{proof}

We are now ready to present one of the main results of  this paper, which characterizes the existence of the spray scalar $H$.
\begin{theorem}\label{Th:dually flat}
Let $ (M,F) $ be a Finsler manifold with geodesic spray coefficients  $G^i$. Then,  the $S$-scalar $H$ exists if and only if $(M,F)$ is  dually flat.
\end{theorem}

\begin{proof}
Assume that the $S$-scalar $H$ exists, then $H=\frac{1}{12} y^k \partial_k F^2$, by Proposition \ref{Prop:S-scalar}. By differentiating  with respect to $y^i$, we obtain 
$$ 12 H_i=y^r\partial_r \dot{\partial}_iF^2+\partial_iF^2.$$
Plugging $H_i$ from \eqref{Eq:H_i} into the above equation, we get
$$ y^r\partial_r \dot{\partial}_iF^2-2\partial_iF^2=0.$$
Moreover, we can write
$$G^i=g^{ij}H_j=-\frac{1}{2}g^{ij}\dot{\partial}_jK,$$
where we have set $K=-2H$. Therefore, \eqref{eq:Dually_flat_1} and \eqref{eq:Dually_flat_2} are satisfied and $(M ,F)$ is dually flat.

Conversely, let $(M ,F)$ be  dually flat, then \eqref{eq:Dually_flat_1} and \eqref{eq:Dually_flat_2} hold. By \eqref{eq:Dually_flat_2}, there exists a function $K\in C^\infty(\mathcal{T} M )$, $3$-homogeneous in $y$, such that 
$$
G^i=-\frac{1}{2} g^{i j}\dot{\partial}_j K.
$$
Then, we have by \eqref{Eq:H_i} 
$$y^r\dot{\partial}_i\partial_r F^2-\partial_i F^2=-2\,\dot{\partial}_i K.$$
In view of \eqref{eq:Dually_flat_1}, the above equation reduces to 
$$\partial_iF^2=-2\,\dot{\partial}_iK.$$
Contracting the last equation by $y^i$, we have
$$K=-\frac{1}{6}y^r\partial_rF^2.$$
That is,  the $S$-scalar $H$ exists and is given in terms of $K$ by $H=-\frac{1}{2} K=\frac{1}{12}y^r\partial_rF^2$.
\end{proof}

It should be noted that, due to  the above theorem, the geometry of     dually flat metrics can be built from a certain scalar function defined on the tangent bundle, the spray scalar $H$.

\begin{proposition}
Let $(M,F)$ be a Berwald manifold. Then, the $F$-covarient coefficients   $H_i$ are quadratic in $y$ if and only if $p^r_{jk}(x)g_{rh}(x,y)$ are functions of $x$ only for some functions $p^i_{jk}(x)$.
\end{proposition}

\begin{proof}
Let $(M,F)$ be   Berwaldian, then $S$ is quadratic and hence we can write $G^i=p^i_{jk}(x)y^jy^k$. 
Using the definition of $H_i$, we have
$$H_i=g_{ir}G^r=g_{ir}p^r_{jk}(x)y^jy^k.$$
So, $H_i$ is quadratic if and only if $g_{ir}p^r_{jk}$ are functions of $x$ only. 
\end{proof}

\begin{remark} 
\label{Remark_semi_field}
It should be noted that the property that  $g_{ir}p^r_{jk}$ are functions of $x$   implies that  $p^r_{jk}C_{rih}=0$. There are examples of Finsler metrics satisfying   the condition $p^r_{jk}C_{rih}=0$ (see \cite{semi-concurrent}).  
\end{remark}
As an example of a Finsler metric whose $F$-covarient coefficients   $H_i$ are quadratic, we have:
\begin{example}  
Let $(M,F)$ be a Finsler manifold where $M=\mathbb{R}^4$ and $F$ is given by
$$F=\sqrt{\sqrt{y_1^4+y_2^4+y_3^4}+x_4y_4^2}.$$
The only non-zero geodesic  spray coefficient is given by
$$G^4=\frac{1}{4} \frac{y_4^2}{x_4}.$$
Also, the only non-zero $F$-covarient coefficient is given by
$$H_4=\frac{1}{4}y_4^2.$$
Hence, $H_i$ are quadratic in $y$.
Here, we have
$$ g_{i4}p^i_{44}= g_{44}p^4_{44}=x_4 \left( \frac{1}{4} \frac{1}{x_4}\right)=\frac{1}{4}.$$
\end{example}

\begin{lemma} Let $(M,F)$ be a Finsler manifold and let $H_i$ be the $F$-covariant coefficients of the geodesic spray of $F$.  The   objects $H_{ij}:=\dot{\partial}_jH_i$ and $H_{ijk}:=\dot{\partial}_kH_{ij}$ have the following formulas in terms of the Finsler function $F$
\begin{eqnarray*}
  H_{ij}&=& \frac{1}{4}(2y^r\partial_rg_{ij}+ \partial_j\dot{\partial}_iF^2- \partial_i\dot{\partial}_jF^2) \\
  H_{ijk}&=& \frac{1}{2}(\partial_kg_{ij}+\partial_jg_{ik}-\partial_ig_{jk})+y^r\partial_rC_{ijk}.
\end{eqnarray*}
\end{lemma}

The following Proposition gives some properties for $H^i_{jk}:=H_{rjh} \ g^{ri}$.
\begin{proposition} \label{cor.1} The $n^3$ functions $H^i_{jk}(x,y):=H_{rjh} \ g^{ri}$ have the following  properties
\begin{description}
   \item[(a)] $H^i_{jk}=\gamma^i_{jk}+g^{ir}y^h\partial_hC_{rjk}$, where ${\gamma}^i_{jk}=\frac{1}{2}{g}^{ir}(\partial_j{g}_{kr}
   +\partial_k{g}_{rj}-\partial_r{g}_{jk})$.
  \item[(b)] $H^i_{jk}y^k=G^i_{j}+2G^r C^i_{rj}$.
  \item[(c)] $H^i_{jk}y^j y^k=2G^i$.
\end{description}
\end{proposition}

\begin{remark} In view of Proposition \ref{cor.1}\textbf{(c)}, the $n^3$ functions $H^i_{jk}(x,y):=H_{rjk}g^{ri}$ have some  properties  similar to those of  Berwald and Cartan connections. However,  by applying a change of coordinates,  $H^i_{jk}$ do not obey the transformation law of a linear connection.  
\end{remark}
Now, we define the following  objects:
\begin{eqnarray*}
  H_{ijkh} &:=&\dot{\partial}_h H_{ijk},\\
  \mathcal{L}_{jkh}&:=&y^iH_{ijkh}.
\end{eqnarray*}

\begin{proposition} \label{prop.1} Let $(M,F)$ be a Finsler manifold. Let $G^h$ and $G^h_i:=\dot{\partial}_i G^h$  be respectively the geodesic spray and the Barthel connection attached to  $F$. The   object $H_{ijkh}(x,y)$, defined above, has the following properties
\begin{description}
   \item[(a)]  $H_{ijkh}$ is symmetric with respect to $j,k$ and $h$.
   \item[(b)]  $H_{ijkh}= \partial_jC_{ikh}+\partial_kC_{i hj}+\partial_hC_{ijk}-\partial_iC_{jkh}+y^r\partial_rC_{ijkh}.$
   \item[(c)]  $H_{ijkh} \, y^j=y^j \partial_jC_{ikh}$ and $ H_{ijkh} \, y^i=\mathcal{L}_{jkh}=-y^i \partial_iC_{jkh}.$
   \item[(d)]  $H_{ijkh}= 2C_{irjkh} G^r + 2\mathfrak{S}_{jkh}\set{C_{irjk}G^r_h+C_{irj}G^r_{kh}}+G_{jikh}$.
     \item[(e)]  $H_{ijkh}\,y^i= -4C_{rjkh} G^r -4\mathfrak{S}_{jkh}\set{C_{rjk}G^r_h}+G_{jikh}\,y^i$,
\end{description}
where, $C_{ijkh}:=\dot{\partial}_h C_{ijk}$,   $C_{ijkhr}:=\dot{\partial}_r C_{ijkh}$, $G^h_{ij}:=\dot{\partial}_j G^h_i$ is the Berwald connection and     $G_{jikh}$ is the hv-Berwald curvature tensor.
\end{proposition}

\begin{remark} In view of the above Proposition, the   object $H_{ijhk}$ satisfies properties similar to those satisfied by the hv-Berwald curvature tensor $G_{ijkh}$. Unfortunately, $H_{ijkh}$ are not  tensor components. So, we will call $H_{ijkh}(x,y)$ and $\mathcal{L}_{jkh}$ Berwald and Landsberg functions, respectively.
\end{remark}

Now, let us define the following new special Finsler spaces.
\begin{definition} Let $(M,F)$ be a Finsler manifold and $H_i$ the associated $F$-covariant coefficients. Then, $(M,F)$ is called $H$-Berwald  space if $H_{ijkh}$ vanishes, that is, if $H_i$  is quadratic. Also, $(M,F)$ is called  $H$-Landsberg  space if $\mathcal{L}_{jkh}$ vanishes.
\end{definition}

\noindent\textbf{Convention:} The following Relations
\begin{eqnarray*}
  C_{irjkh} G^r &+& \mathfrak{S}_{jkh}\set{C_{irjk}G^r_h+C_{irj}G^r_{kh}}=0 \\
 C_{rjkh} G^r & +& \mathfrak{S}_{jkh}\set{C_{rjk}G^r_h}=0
\end{eqnarray*}
will be called I-condition and II-condition, respectively.

\vspace{5pt}
We have the following   theorem.
\begin{theorem} For a Finsler manifold $(M,F)$, we have
 \begin{description}
    \item[(a)] Every  $H$-Berwald   space is an $H$-Landsberg   space, and the converse is not true in general.
    \item[(b)] Under the I-condition, the two notions of   $H$-Berwald    and   Berwald spaces coincide. 
        \item[(c)] Under the II-condition, the two notions of   $H$-Landsberg     and    Landsberg  spaces coincide. 
\end{description}
\end{theorem}
\begin{proof} The proof is straightforward and therefore omitted.
\end{proof}

 Considering  the above Theorem, we have the  diagram:

{{{\begin{displaymath}
    \xymatrix{
\text{H-Berwald    space}
 \ar@{=>}[dd]_{(a)}\ar@{<=>}[rr]^{I-condition}_{(b)}& &\text{Berwald  space} \ar@{=>}[dd]\\&&\\
\text{H-Landsberg   space } \ar@{<=>}[rr]^{II-condition}_{(c)}&& \text{Landsberg  space}}
\end{displaymath}}}}

\begin{theorem} \label{Th:H-Landsberg}
Let $F=F(x, y)$  be a Finsler metric.  If the $S$-scalar $H$ exists, then $F$    is  an $H$-Landsberg metric.
\end{theorem} 

\begin{proof}
As $H$ exists, then $H_i=\dot{\partial}_i H$. Therefore, by the homogeneity of $H$, we have
$$\mathcal{L}_{jkh}=y^iH_{ijkh}=0.$$
That is, $F$ is  an $H$-Landsberg metric.
\end{proof}

In view of the above Theorem \ref{Th:H-Landsberg} and Theorem \ref{Th:dually flat}, we have
\begin{corollary}\label{Th:H-Landsberg & dually flat}
Let $F=F(x, y)$  be a Finsler metric.  If the $S$-scalar $H$ exists, then $F$ is both $H$-Landsbergian and dually flat metric.
\end{corollary}

In contrast to the unicorn Landsberg problem, some examples of H-Landsberg Finsler metrics which are not H-Berwaldian will be given in Section 5.

\section{Special Finsler spaces}
In this section, we study two special Finsler metrics, namely, projectively flat metrics and spherically symmetric Finsler metrics, in relation to the $F$-covariant coefficients $H_i$. Let's start with the projectively flat metrics.

\subsection{Projectively flat Finsler metrics}
For the  projectively flat  metrics, we study the case when $H_i=P(x,y)y_i$ provided that $P(x,y)$ is positively homogeneous function of degree $1$ in $y$. 

\begin{definition}\cite{Chern-Shen-book}
 A Finsler metric $F$ on a manifold $M$ is said to be  projectively flat if and only if $F$ satisfies the property
\begin{equation}
\label{Eq:Projectivly-flat}
y^k\dot{\partial}_i\partial_k F-\partial_i F=0 .
\end{equation}
In this case, $G^i=P(x, y) y^i $ with $P=\frac{y^i \partial_i F }{2 F}$ and $P$ is called  the projective factor of $F$.
\end{definition}

\begin{proposition} A Finsler metric $F$ on $M$ is projectively flat (i.e., $G^i=P(x, y) y^i$) if and only if 
the $F$-covariant coefficients $H_i$ of the geodesic spray $S$ are given by $H_i =P(x,y) y_i \, ; \   y_i:=g_{ij}\,y^j$.
\end{proposition}

\begin{lemma} Consider a    projectively flat Finsler metric $F$  with $F$-covariant coefficients $H_i=P(x,y)y_i$  and   projective factor $P(x,y)$.  Then, we have
\begin{eqnarray*}
  H_{ij}&=& P_j \, y_i+P\, g_{ij} ,\\
  H_{ijk}&=& P_{jk} y_i + P_j \, g_{ik}+ P_k \, g_{ij}+2P \, C_{ijk},
\end{eqnarray*}
where $P_j:=\dot{\partial}_j P $ and $P_{jk}:=\dot{\partial}_k P_j$.
\end{lemma}

\begin{corollary}  \label{Cor4.4}  In the case of projectively flat metric,  $H^i_{jk}(x,y):=H_{rjh} \ g^{ri}$  satisfies
\begin{description}
   \item[(a)] $H^i_{jk}=P_{jk} y^i+P_j \delta^i_k+P_k \delta^i_j+2P \, C^i_{jk}$
  \item[(b)] $H^i_{jk}y^k=2P\delta^i_j$,\quad    $H^i_{jk}y^j y^k= 2P\, y^i$.
  \end{description}
\end{corollary}

\begin{corollary} For a projectively flat metric,  in view of Corollary \ref{Cor4.4} (a), we have
$$H^i_{jk}=G^i_{jk}+2P \, C^i_{jk},$$
where $G^i_{jk}=P_{jk} y^i+P_j \delta^i_k+P_k \delta^i_j$. In this case  $H^i_{jk}$ is   a linear connection.
\end{corollary}

\begin{corollary} For the projectively flat metrics, the Berwald functions $H_{ijkh}(x,y)$ has the following properties
\begin{description}
        \item[(a)]  $H_{ijkh}= 2P_{jkh} y_i+ 2\mathfrak{S}_{jkh}\set{P_{jk}\,g_{ih}+P_j \,C_{irj}}+2P \, C_{ijkh}$.
        \item[(b)]  $H_{ijkh}= G_{jikh}+ 2\mathfrak{S}_{jkh}\set{P_j \,C_{ikh}}+2P \, C_{ijkh}$.
     \item[(c)]  $H_{ijkh}\,y^i= G_{jikh} y^i -4P\, C_{jkh}$.
\end{description}
\end{corollary}

\begin{remark} As we have seen before, the Berwald functions $H_{ijkh}$ are not components of a tensor field. However, in the case of projectively flat metrics, they are components of a tensor field of type $(0,4)$.
\end{remark}

\subsection{Spherically symmetric metrics}

~\par In this subsection, we investigate the spherically symmetric Finsler metrics $F$ in relation to the $F$-covariant coefficients $H_i$ of $F$.

A spherically symmetric Finsler   $F$ on $\mathbb{B}^n\left(r_0\right) \subset \mathbb{R}^n$ is defined as follows (see \cite{Guo-Mo}):
$$F(x,y)=|y|\ \varphi\left(|x|,\frac{\langle x, y \rangle}{|y|}\right),$$
where $\varphi:[0,r_0)\times\mathbb{R}^n\to \mathbb{R}$ is a smooth function on $\mathcal{T}M$, $(x,y)\in T\mathbb{B}^n(r_0)\backslash \{0\}$ and $|\cdot|$, $\langle \cdot , \cdot \rangle$ are the standard Euclidean norm and inner product on $\mathbb{R}^n$, respectively.

 $F$ can   simply be written in the form  $$F=u\ \varphi(r,s)$$
  where $u=|y|$, $r=|x|$ and $s=\frac{\langle x, y \rangle}{|y|}$.

\medskip
Also, we have the following notations
$$y_i:=\delta_{ih} y^h, \quad x_i:=\delta_{ih} x^h.$$
That is, $y_i$ and $x_i$ are equivalent to $y^i$ and $x^i$, respectively. Therefore, we conclude that $y_i \neq F \frac{\partial F}{\partial y^i}$, but rather $y_i=u \frac{\partial u}{\partial y^i}$. Moreover, we have  
$$y^iy_i=u^2, \quad x^ix_i=r^2, \quad y^ix_i=x^iy_i=\langle x,y \rangle.$$

 \medskip

The components $g_{ij}$ of the metric tensor attached to      $F=u \varphi(r,s)$ are given by  
\begin{equation}
\label{Eq:g^ij}
g_{ij}=\sigma_0\ \delta_{ij}+\sigma_1\  x_ix_j+\frac{\sigma_2}{u} (x_iy_j+x_jy_i)+\frac{\sigma_3}{u^2}y_iy_j,
\end{equation}
where the functions $\sigma$'s are defined by  $$\sigma_0:=\varphi(\varphi-s\varphi_s),\quad \sigma_1:= \varphi_s^2+\varphi\varphi_{ss},\quad \sigma_2:= (\varphi-s\varphi_s)\varphi_s-s\varphi\varphi_{ss},  \quad \sigma_3:= s^2\varphi\varphi_{ss}-s(\varphi-s\varphi_s)\varphi_s.$$
Clearly, $$s \sigma_2+\sigma_3=0, \quad s \sigma_1+\sigma_2=\varphi \varphi_s, $$
 where the subscript $s$ denotes the derivative with respect to $s$. Spherically symmetric Finsler metrics have been extensively studied by various authors. For further details, we refer to \cite{Guo-Mo, Zhou_Mo}.

  The coefficients  $G^i$ of the geodesic spray  of $F$ are  given by
\begin{equation}\label{G}
  G^i=uPy^i+u^2 Qx^i,
\end{equation}
where the functions $P$ and $Q$ are defined by
\begin{equation}\label{P,Q}
 Q:=\frac{1}{2 r}\frac{ -\varphi_r+s\varphi_{rs}+r\varphi_{ss}}{\varphi-s\varphi_s+(r^2-s^2)\varphi_{ss}}, \quad P:=-\frac{Q}{\varphi}(s\varphi +(r^2-s^2)\varphi_s)+\frac{1}{2r\varphi}(s\varphi_r+r\varphi_s).
\end{equation}

Now, the $F$-covariant coefficients $H_i$ of the Finsler metric $F=u\varphi(r,s)$ can be calculated, using \eqref{G} and \eqref{P,Q},  as follows:
\begin{align*}
H_j&=g_{ij}G^i\\
&=(\sigma_0\ \delta_{ij}+\sigma_1\  x_ix_j+\frac{\sigma_2}{u} (x_iy_j+x_jy_i)+\frac{\sigma_3}{u^2}y_iy_j)(uPy^i+u^2 Qx^i)\\
&=(uP(\sigma_0+s \sigma_2+\sigma_3)+uQ(r^2\sigma_2+s\sigma_3))y_j\\
&+(u^2P(s\sigma_1+\sigma_2)+u^2Q(\sigma_0+r^2\sigma_1+s\sigma_2))x_j.
\end{align*}
Therefore, we have
$$H_j=u(P\sigma_0+(r^2-s^2)\sigma_2Q)y_j+u^2(P\varphi\varphi_s+(\varphi^2+(r^2-s^2)\sigma_1)Q)x_j.$$

A spherically symmetric   metric $F=u\varphi(r,s)$ is both projectively flat and dually flat if and only if the function $\varphi$ is given by \cite{Najafi}:
$$
\varphi(r,s)=\frac{\sqrt{\left(\mathrm{k}^2-c^2 r^2\right)+c^2 s^2}+c s}{\left(\mathrm{k}^2-c^2 r^2\right)},
$$
where $k$ and $c$ are constants.\\
In this case, we have 
\begin{align}\label{not quadratic}
H_j&=\frac{1}{2} \frac{u c\left(\sqrt{-c^2 r^2+c^2 s^2+k^2}+c s\right)^2}{\left(c^2 r^2-k^2\right)^2 \sqrt{-c^2 r^2+c^2 s^2+k^2}} y_j-\frac{1}{2} \frac{\left(\sqrt{-c^2 r^2+c^2 s^2+k^2}+c s\right)^3 c^2 u^2}{\left(c^2 r^2-k^2\right)^3\sqrt{-c^2 r^2+c^2 s^2+k^2}} x_j
\end{align}
In view of the fact that $F$ is dually flat, the spray scalar $H$ exists, by Theorem \ref{Th:dually flat}, and can be obtained by the following computation.\pagebreak
\begin{align*}
y^j H_j&=\frac{1}{2} \frac{u c\left(\sqrt{-c^2 r^2+c^2 s^2+k^2}+c s\right)^2}{\left(c^2 r^2-k^2\right)^2 \sqrt{-c^2 r^2+c^2 s^2+k^2}} u^2-\frac{1}{2} \frac{\left(\sqrt{-c^2 r^2+c^2 s^2+k^2}+c s\right)^3 c^2 u^2}{\left(c^2 r^2-k^2\right)^3\sqrt{-c^2 r^2+c^2 s^2+k^2}} su\\
&=-\frac{1}{2}\frac{c u^3 \left(\sqrt{-c^2 r^2+c^2 s^2+k^2}+c s\right)^3}{\left(c^2 r^2-k^2\right)^3}\\
&=3H.
\end{align*}
That is, $H=-\frac{1}{6}\frac{c u^3 \left(\sqrt{-c^2 r^2+c^2 s^2+k^2}+c s\right)^3}{\left(c^2 r^2-k^2\right)^3}$.\\

It should be noted that dually flat Finsler metrics (for which the $S$-scalar exists) are $H$-Landsberg metrics, by Theorem \ref{Th:H-Landsberg}. On the other hand, projectively flat and dually flat spherically symmetric Finsler metrics are not $H$-Berwaldian since the $F$-covarient coefficients $H_i$ of such spaces are not quadratic in $y$ as shown in Equation \eqref{not quadratic}.

Consequently, we have the following interesting result.
\begin{theorem}
The class of projectively flat and dually flat spherically symmetric Finsler metrics constitutes a solution for the $H$-unicorn Landsberg problem.
\end{theorem}

\section{Examples}

In this section, we provide two  examples. The first one is   an $n$-dimensional projectively flat $H$-Landsberg Finsler metric which is  non $H$-Berwaldian. The second one is an $H$-Berwald metric.

\begin{example}
Let $M=\mathbb{R}^n$   and $F$ be a Randers metric  given by
$$F(x,y)=\alpha+\beta= c|y|+\langle a+x,y\rangle,$$
where the Riemannian metric $\alpha=c|y|$, the one form $\beta=\langle a+x,y\rangle$    ($c\in \mathbb{R}$   and  $a$ is a fixed vector in $\mathbb{R}^n$) and $|\cdot|$,   $\langle \cdot , \cdot \rangle$ are the standard Euclidean norm and inner product on $\mathbb{R}^n$, respectively.  Keeping in mind that $x_i$ and $y_i$ are the same as $x^i$ and $y^i$, respectively, since the background metric is the Kronecker delta, then  we can rewrite $F$ as follows:
$$F(x,y)=c\sqrt{y_1^2+y_2^2+\cdots + y_n^2}+(a_1+x_1)y^1+(a_2+x_2)y^2+\cdots+(a_n+x_n)y^n,$$

The Randers metric $F$ is a projectively flat metric. Indeed, we have
$$\partial_j F=y_j, \quad \dot{\partial}_k\partial_jF=\delta_{jk}.$$
Also
$$\partial_j F^2=2F y_j, \quad \dot{\partial}_k\partial_jF^2=2\ell_k y_j+2F\delta_{jk}.$$
Then,   substituting into \eqref{Eq:Projectivly-flat},  we get
$$y^j \dot{\partial}_k\partial_jF-\partial_k F=y^j \delta_{jk}-y_k=0.$$
That is,   $F$ is projectively flat. Moreover, the $F$-covariant coefficients $H_i$ are given by
$$H_i= \frac{1}{4}(y^r\partial_r\dot{\partial}_iF^2 - \partial_iF^2)=\frac{1}{2} |y|^2\ell_i.$$
But $\ell_i=\frac{c}{|y|}y_i+a_i+x_i$. Then, $H_i$ can be written as  
$$H_i= \frac{c}{2}\,|y| y_i+\frac{1}{2} |y|^2(a_i+x_i).$$
Taking the derivative with respect to $y^j$, we obtain
$$H_{ij}=\frac{c}{2|y|}y_iy_j+\frac{c}{2}|y|\delta_{ij}+y_j(a_i+x_i).$$
Differentiating $H_{ij}$ with respect to $y^k$,  we get
$$H_{ijk}=-\frac{c}{2|y|^3}y_iy_jy_k+\frac{c}{2}\,|y|(y_i\delta_{jk}+y_j\delta_{ik}+y_k\delta_{ij})+\delta_{jk}(a_i+x_i).$$
Similarly, the differentiation   with respect to $y^h$, the $H$-Berwald tensor can be written in the following form
\begin{align*}
H_{ijkh}&=\frac{3c}{2|y|^5}y_iy_jy_ky_h-\frac{c}{2|y|^3}(y_iy_j\delta_{kh}+y_iy_k\delta_{jh}+y_jy_k\delta_{ih}\\
&+y_iy_h\delta_{jk}+y_jy_h\delta_{ik}+y_ky_h\delta_{ij})+\frac{c}{2|y|}(\delta_{ih}\delta_{jk}+\delta_{jh}\delta_{ik}+\delta_{kh}\delta_{ij})
\end{align*}
Now, the $H$-Landsberg tensor is given by$$\mathcal{L}_{jkh}=y^iH_{ijkh}=0.$$
Therefore, the metric $F$ is $H$-Landsbergian. On the other hand, $F$ is non $H$-Berwaldian because $H_i$ are not quadratic in $y$.
\end{example}

Finally, we provide the following  example of a class of $H$-Berwald spaces.
 
\begin{example}
Let $M=\mathbb{R}^3$ and $F$ be defined as follows
$$F(x,y)=\sqrt{\frac{a_1y_1^4+a_2y_1^2y_3^2+a_3y_2^2y_3^2}{b_1y_1^2+b_2y_2^2+b_3y_3^2}+f_1(x_3)y_1^2+f_2(x_3)y_2^2}.$$
Using \eqref{Eq:H_i}, $F$-covariant coefficients $H_i$ are given by
$$H_1=2y_1y_3\frac{d f_1(x^3)}{dx^3}, \quad H_2=2y_2y_3\frac{d f_2(x^3)}{dx^3}, \quad H_3= -y_1^2\frac{d f_1(x^3)}{dx^3}-y_2^2\frac{d f_2(x^3)}{dx^3}.$$
Since the coefficients $H_i$ are quadratic in $y$, the metric $F$ is $H$-Berwaldain.
\end{example}

\section*{Conclusion}
We conclude this work with the following remarks:

$\bullet$  A spray manifold $(M, S)$ is a manifold $M$ equipped with a spray $S$. For a spray manifold $(M, S)$, by starting with a spray $S$ defined by its coefficients $G^i$, we can construct a geometry on the manifold. We can obtain the canonical nonlinear connection $G_j^i=\paa_j G^i$. This nonlinear connection allows us to define the horizontal sub-bundle $H T M$. Furthermore, we can investigate the induced Berwald connection and its associated torsion and curvature tensors.

\medskip

$\bullet$ 
For a Finsler   manifold $(M, F)$, we have investigated    and uncovered additional geometric structures through the concept of $F$-covariant coefficients $H_i$ of the geodesic spray $S$ of $F$. For instance, we have introduced the notions of H-Berwald and H-Landsberg spaces. Analogous to Berwald spaces, where the coefficients $G^i$ are quadratic in $y, \mathrm{H}$-Berwald spaces are characterized by quadratic $H_i$ coefficients. In contrast to the longstanding open problem of finding non-trivial Landsberg metrics that are not Berwaldian, we have provided examples of H-Landsberg Finsler metrics that are not H-Berwaldian.

\medskip

$\bullet$ As a continuation of our investigation of the $F$-covariant coefficients, analogous to the Berwald connection $G_{i j}^h$ obtained from the coefficients $G^i$, in the case of projectively flat metrics  the functions $H_{i j}^h$ obtained from the $F$-covariant coefficients $H_i$ define a linear connection.

\section*{Declarations} 
\begin{itemize}
\item \textbf{Competing interests}: The authors  declare  no conflict of interest.
  \item \textbf{Availability of data and material}: Not applicable.
  \item \textbf{Funding}: Not applicable.
  \item \textbf{Authors' contributions}: The authors  have made substantive contributions
to the article and assume full responsibility for its content. The authors  read and approved the final manuscript.

\end{itemize}


\begin{thebibliography}{99}



\bibitem{sprays}
W. Ambrose, R. S. Palais and I. M. Singer,  \emph{Sprays}, Acad. Brasil. Ciencas,  \textbf{32}   (1960),  163-178.



\bibitem{Chern-Shen-book}
 S. Chern and Z. Shen,  \emph{Riemann-finsler geometry},  World Scientific, Singapore, 2005.



\bibitem{Elgendi-SSM}
 S. G. Elgendi, \textit{ On the classification of Landsberg spherically symmetric
Finsler metrics, } Int. J. Geom. Methods Mod. Phys., \textbf{18}, 14 (2021) 2150232. arXiv:  2110.07252 [math.DG]. 



\bibitem{Elgendi-SSM-2}
 S. G. Elgendi, \textit{A note on “On the classification of Landsberg spherically symmetric Finsler metrics” } Int. J. Geom. Methods Mod. Phys., \textbf{20}, 6 (2023) 2350096-300.        
arXiv:  2302.09848 [math.DG].


\bibitem{r21}
J.~Grifone, \emph{Structure pr\'esque-tangente et connexions,
\textsc{I}}, Ann.
  Inst. Fourier, Grenoble,\textbf{ 22}, 1 (1972), 287-334.


\bibitem{Guo-Mo}
E. Guo and X. Mo, \emph{The geometry of spherically symmetric Finsler manifolds}, Springer, 2018.


\bibitem{Zhou_Mo}
X. Mo and L. Zhou, \emph{The curvatures of spherically symmetric Finsler metrics in $\mathbb{R}^n$}, (2014). arXiv: 1202.4543v4 [math.DG].


\bibitem{Najafi}
B. Najafi,  \emph{On projectively flat spherically symmetric Finsler metrics}, (2015).  arXiv: 1503.05415 [math.DG].  


\bibitem{Shen-dually}
Z. Shen, \emph{Riemann-Finsler geometry with applications to information geometry},  Chinese Annal.  Math.,  Series B, \textbf{27} (2006), 73-94.

 

\bibitem{Szilasi-book} J. Szilasi, R.L. Lovas, D.Cs. Kert\'esz:
  \emph{Connections, Sprays and Finsler Structures}, World Scientific, 2014.

\bibitem{NF_Package}
Nabil~L. Youssef and S. G. Elgendi, \emph{New Finsler   package},  Comput. Phys. Commun., \textbf{185}, 3 (2014), 986--997. arXiv:  1306.0875 [math.DG]

\bibitem{semi-concurrent}
Nabil L. Youssef, S. G. Elgendi and Ebtsam H. Taha, \emph{Semi-concurrent vector field in Finsler geometry},  Differ. Geom. Appl.,  \textbf{65} (2019), 1-15.   arXiv: 1802.02405 [math.DG]. 



  \end{thebibliography}
\end{document}